\documentclass{proc-l}

\copyrightinfo{2006}{American Mathematical Society}

\newtheorem{lemma}{Lemma}
\newtheorem{theorem}{Theorem}
\newtheorem{corollary}{Corollary}

\numberwithin{equation}{section}

\begin{document}


\title[Asymptotic expansion for inverse moments of distributions]{Asymptotic expansion for inverse moments of binomial and Poisson distributions}

\author{Marko \v Znidari\v c}
\address{Department of Quantum Physics, University of Ulm, Germany {\em and}\\
Faculty of Mathematics and Physics, University of Ljubljana, Slovenia}
\curraddr{}
\email{}
\thanks{}


\date{}

\dedicatory{}


\begin{abstract}
An asymptotic expansion for inverse moments of positive
binomial and Poisson distributions is derived. The expansion coefficients of the asymptotic series are
given by the positive central moments of the distribution. Compared to
previous results, a single expansion formula covers all (also non-integer) inverse moments. In addition, the approach can be generalized to other positive distributions.
\end{abstract}

\maketitle

\section{Introduction}

Inverse moments of probability distributions can arise in several contexts. In particular, they are relevant in various statistical applications, see {\em e.~g.} Grab and Savage~\cite{Grab:54}, Mendenhall and Lehman~\cite{Mendenhall:60}, Jones and Zhigljavsky~\cite{Jones:04}, and references therein. Recently it has been shown~\cite{Znidaric:05} that the first two inverse moments of positive binomial distribution are needed for the calculation of the running time of a particular quantum adiabatic algorithm solving 3-SAT problem. Because closed expressions are usually not possible it is of interest to derive asymptotic expansions.

In the present paper we are going to study the inverse moments of positive binomial and Poisson distributions. For binomial distribution,
\begin{equation}
P(x=i)=\frac{1}{1-q^n} {n \choose i} p^i q^{n-i},\qquad i=1,\ldots,n
\end{equation}
with $q=1-p$, the $r$-th inverse moment is given by
\begin{equation}
f_r(n)=(1-q^n)\mathbb{E}(1/x^r)=\sum_{i=1}^n{{n \choose i} p^i q^{n-i} \frac{1}{i^r}}.
\label{eq:def}
\end{equation}
Similarly, for Poisson distribution the $r$-th moment $g_r$ is
\begin{equation}
g_r=\sum_{s=1}^\infty{{\rm e}^{-m} \frac{m^s}{s!} \frac{1}{s^r}}.
\label{eq:defP}
\end{equation}
We are going to derive an asymptotic expansion of $f_r(n)$ and $g_r$ in terms of inverse powers of $(np+q)$ or $(m+1)$ for binomial and Poisson distribution, respectively. The expansions work for arbitrary real and positive $r$. We will also generalize the result to arbitrary positive distributions. The expansion coefficients are given by the central moments of the distribution in question which is intuitively appealing result. Namely, for sufficiently sharp distribution the inverse moment will be approximately given by the inverse of the average value, {\em i.e.} $\mathbb{E}(1/x^r) \approx 1/\mathbb{E}(x)^r$. If the distribution is non-zero at a single point, then this result is exact. Higher orders in the asymptotic expansion should therefore ``measure'' how much the distribution in question differs from a point distribution, $P(x=i)=\delta_{i,\mathbb{E}(x)}$. This is the reason why it is natural that the expansion is given in terms of central moments. But let us first briefly review known results about inverse moments.     

The problem of calculating the inverse moments of positive distributions, {\em e.~g.} $f_r(n)$ or $g_r$, has a rather long history. Already in 1945 Stephan~\cite{Stephan:45} studied the first and second inverse moment for binomial distribution. Grab and Savage~\cite{Grab:54} calculated tables of reciprocals for binomial and Poisson distributions as well as derive a recurrence relation. They also derived an exact expression for the first inverse moment for a Poisson distribution involving an exponential integral. David and Johnson~\cite{David:56} used an approximating formula for binomial distribution in terms of inverse powers of $(np)$. Govindarajulu~\cite{Govi:63} obtained a recursive formula for inverse moments of binomial distribution. The recurrence formula though is rather complicated as it involves recursion in $n$ as well as in the power of the moment $r$. Asymptotic expansions have also been studied. Such expansions are of value also because the formal definitions of $f_r(n)$ (\ref{eq:def}) or $g_r$ (\ref{eq:defP}) are not simple to evaluate for large values of $np$ or $m$. One has to sum many terms, each being a product of a very small and a very large number, all together resulting in a small inverse moment. Tiku~\cite{Tiku:64} derived an asymptotic expansion of $g_r$ using Laguerre polynomials. Chao and Strawderman~\cite{Chao:72} considered slightly different inverse moments defined as $\mathbb{E}(1/(x+a)^r)$ which are frequently easier to calculate. Simple expression for integer $a\ge 1$ and $r=1$ are derived. Kabe~\cite{Kabe:76} derived a general series with the expansion coefficients given by positive factorial moments $\mu^{[i]}=\mathbb{E}(x!/(x-i)!)$ of the distribution. He also derived a general formula
\begin{equation}
\mathbb{E}\left(\frac{1}{(x+a)^r}\right)=\frac{1}{\Gamma(r)} \int_0^\infty{\!\! y^{r-1}{\rm e}^{-a y} M(-y) {\rm d}y},
\label{eq:int}
\end{equation}
where $M(t)$ is a moment generating function, $M(t)=\mathbb{E}({\rm e}^{tx})=\sum_{j=0}^\infty{\mu_j x^j/j!}$ with $\mu_j=\mathbb{E}(x^j)$. For further discussion of negative moments and generating functions see also Cressie {\em et.al.}~\cite{Cressie:81}. Various other identities involving generating functions are derived in Rockower~\cite{Rockower:88}. Equation~\ref{eq:int} will serve as a starting point for our derivation. Quite recently the asymptotic expansions of the first inverse moment for binomial distribution have been considered by Marciniak and Weso\l owski~\cite{Marciniak:99} and by Rempala~\cite{Rempala:03}. In~\cite{Marciniak:99} an expansion in terms of Euler polynomials is given while in~\cite{Rempala:03} a more elegant expansion in terms of inverses of $n^{[j]}=n!/(n-j)!$ is presented. In both works only the expansion of the first moment is given. Jones and Zhigljavsky~\cite{Jones:04} gave a general method using Stirling numbers to derive the asymptotic expansion for Poisson distribution for an arbitrary real power $r$. The expansion though gets more complicated for larger $r$.

In the present work we are going to systematically derive a general asymptotic expansions of $f_r(n)$ and $g_r$ for an arbitrary real power $r$. What is more, the expansion for an arbitrary $r$ will be given by a single simple formula. The paper is organized as follows. In section~\ref{sec1} we derive asymptotic expansion for binomial distribution. In section~\ref{sec2} a similar result is derived for Poisson distribution. Finally, in section~\ref{sec3} the method is generalized to arbitrary positive distributions.
  
\section{Asymptotic Expansion for Binomial Distribution}
\label{sec1}

Note that for integer $r$ the inverse moment $f_r(n)$ can be formally written in terms of a generalized hypergeometric function~\cite{Abramowitz:72}, 
\begin{eqnarray}
f_1(n)&=& np q^{n-1}\> _2F_3\left( {1,1,1-n \atop 2,2};-\frac{p}{q}\right) \nonumber \\
f_2(n)&=& np q^{n-1}\> _3F_4\left( {1,1,1,1-n \atop 2,2,2};-\frac{p}{q}\right),
\label{eq:gamma_F}
\end{eqnarray}
and analogously for larger $r$'s. Although exact, these expressions are not very illuminating. In fact, using a series expansion, {\em e.g.} 
\begin{equation}
_2F_3 \left( {a,b,c \atop d,g};z\right)=\sum_{k=0}^\infty{\frac{(a)_k (b)_k (c)_k}{ (d)_k (g)_k}\frac{z^k}{k!}},
\end{equation}
where $(a)_k=\Gamma(a+k)/\Gamma(a)$ is the Pochhammer symbol, one can see that  the corresponding series for a generalized hypergeometric function in equation~\ref{eq:gamma_F} is nothing but the sum occurring in the average for $f_r(n)$ (\ref{eq:def}).

Our approach will be different. First we write $f_r(n)$ as an integral. Using identity $\int_0^\infty{x^p {\rm e}^{-r x}{\rm d}x}=\frac{\Gamma(p+1)!}{r^{p+1}}$, with $\Gamma(p+1)=p!$ for integer $p$, and binomial expansion of $(q+p{\rm e}^{-x})^n$ we quickly see that $f_r(n)$ (\ref{eq:def}) is given by
\begin{equation}
f_r(n)=\frac{1}{\Gamma(r)}\int_0^\infty{\!\!\!x^{r-1} \left[(q+p{\rm e}^{-x})^n-1\right]{\rm d}x}.
\end{equation}
This is nothing but a special case of equation~\ref{eq:int}. Integrating once per parts we get
\begin{equation}
f_r(n)=\frac{np}{\Gamma(r+1)}\int_0^\infty{\!\!\!x^r {\rm e}^{-x} (q+p{\rm e}^{-x})^{n-1}{\rm d}x}.
\label{eq:frint}
\end{equation}
The last equation will serve as a starting point for the asymptotic expansion. Before proceeding though, let us prove two auxiliary Lemmas about central moments.

\begin{lemma}
\label{lemma_F}
Let $M(t)$ be generating function of probability distribution,
\begin{equation}
M(t)=\sum_{k=0}^\infty{\mu_k' \frac{t^k}{k!}},
\end{equation}
where $\mu_k':=\mathbb{E}(x^k)$ is $k$-th moment about origin. Then the generating function $F(t)$ of central moments, {\em i.e.} 
\begin{equation}
F(t)=\sum_{k=0}^\infty{ \mu_k \frac{t^k}{k!}},\qquad \mu_k:=\mathbb{E}((x-\bar{x})^k),
\end{equation}
with $\bar{x}=\mathbb{E}(x)$ being the mean, is given by
\begin{equation}
F(t)={\rm e}^{-t \bar{x}} M(t).
\end{equation}
\end{lemma}
\begin{proof}
Lemma is proved by straightforward expansion,
\begin{equation}
F(t)=\sum_{k=0}^\infty\sum_{l=0}^\infty (-\bar{x})^l \frac{t^l}{l!} \mathbb{E}(x^k) \frac{t^k}{k!}.
\end{equation}
Rewriting summation in terms of $j=l+k$ and $l$ and recognizing binomial expansion of $(x-\bar{x})^j$, we get $F(t)=\sum_{j=0}^\infty \mu_j t^j/j!$. 
\end{proof}
In the special case of binomial distribution Lemma~\ref{lemma_F} gives the following result.
\begin{corollary}
For the binomial distribution generating function $F_n(t)$ of central moments is,
\begin{equation}
F_n(t)=(q {\rm e}^{-tp}+p {\rm e}^{tq})^n,
\label{eq:Fdef}
\end{equation}
where central moments are~\footnote{Note that in the definition of $\mu_k(n)$ the sum runs over $i=0,\ldots,n$, that is including $i=0$, which term is for instance excluded in the definition of $f_r(n)$ (\ref{eq:def}).},
\begin{equation}
\mu_k(n)=\sum_{i=0}^n{{n \choose i} p^i (1-p)^{n-i} (i-np)^k}.
\label{eq:mudef}
\end{equation}
\label{lemmaF}
\end{corollary}
\begin{proof}
Generating function for binomial distribution is $M_n(t)=(q+p{\rm e}^t)^n$ and the mean is $\bar{x}=np$. Using Lemma~\ref{lemma_F} immediately gives the result.
\end{proof}

The first few central moments of binomial distribution are readily calculated and are
\begin{multline}
\mu_0(n)=1,\quad \mu_1(n)=0,\quad \mu_2(n)=np(1-p),\\
\quad \mu_3(n)=np(1-p)(1-2p),\quad \mu_4(n)=np(1-p)(1-6p+6p^2+3np-3np^2).
\end{multline}
Higher moments can be calculated using the following Lemma.
\begin{lemma}
For central moments of binomial distribution $\mu_k(n)$ (\ref{eq:mudef}) the following recursive relation holds:
\begin{equation}
\mu_{k+1}(n)=p(1-p)\left\{\frac{{\rm d}\mu_k(n)}{{\rm d}p}+nk \mu_{k-1}(n) \right\},
\label{eq:murecursive}
\end{equation}
where $\mu_k(n)$ is considered as a function of $p$ only ($q$ is replaced by $1-p$).
\label{lemmarec}
\end{lemma}
\begin{proof}
The relation is derived by differentiating definition of $\mu_k(n)$ (\ref{eq:mudef}) with respect to $p$. Using identity
\begin{equation}
\sum_{i=0}^n{{n \choose i} p^i q^{n-i} i(i-np)^k}=\mu_{k+1}(n)+np\mu_k(n),
\end{equation}
obtained by writing $\mu_{k+1}(n)=\sum_{i=0}^n{{n \choose i} p^i q^{n-i} (i-np)^k (i-np)}$, one obtains the recursive relation (\ref{eq:murecursive}).
\end{proof}

\begin{corollary}
Using Lemma~\ref{lemmarec} we can show that
\begin{enumerate}
\item $\mu_k(n)$ is polynomial of order $k$ in $p$
\item $\mu_k(n)$ is polynomial of order $\lfloor \frac{k}{2} \rfloor$ in $n$
\item Leading order in $n$ of $\mu_{2k}(n)$ is $\mu_{2k}(n)=(2k-1)!!(pqn)^k +\mathcal{O}(n^{k-1})$.
\end{enumerate}
\label{cor123}
\end{corollary}
\begin{proof}
Using induction and recursive relation (\ref{eq:murecursive}), together with the initial condition $\mu_1(n)=0$ and $\mu_2(n)=np(1-p)$, we can see that the order of polynomial $\mu_k(n)$ in $p$ is indeed $k$. Similarly, the order in $n$ of $\mu_{k+1}(n)$ is by $1$ larger than that of $\mu_{k-1}(n)$, therefore $\mu_k(n)\sim \mathcal{O}(n^{\lfloor \frac{k}{2}\rfloor})$. Specifically, one can see that the leading order in $\mu_{2k}(n)$ is $\mu_{2k}(n)\sim (pq)^k (2k-1)!!n^k$.    
\end{proof}

Now we are ready to give the asymptotic expansion of $f_r(n)$ (\ref{eq:def}).
\begin{theorem}
For any real $r>0$ inverse moment $f_r(n)$ (\ref{eq:def}) is given by asymptotic expansion
\begin{equation}
f_r(n)=\frac{np}{(np+q)^{r+1}} \left\{ \sum_{k=0}^{m-1}{(-1)^k
  \frac{\mu_k(n-1)}{(np+q)^k} {r+k \choose r}}+\mathcal{O}(1/n^{\lceil
  \frac{m}{2} \rceil}) \right\}.
\label{eq:theorem}
\end{equation}
The terms $\mu_k(n-1)/(np+q)^k$ are of order $\mathcal{O}(1/n^{\lceil
  \frac{k}{2} \rceil})$, where $\lceil x \rceil$ is the smallest integer
larger than $x$. For non-integer $r$'s the binomial symbol is understood as ${r+k \choose r}=(r+1)(r+2)\cdots(r+k)/k!$.
\label{theoremf}
\end{theorem}

\begin{proof}
Writing $(q+p{\rm e}^{-x})^{n-1}={\rm e}^{-xp(n-1)} (q{\rm e}^{xp}+p{\rm e}^{-xq})^{n-1}$ in integral equation for $f_r(n)$ (\ref{eq:frint}), changing the integration variable to $y=x(np+q)$, we get
\begin{equation}
f_r(n)=\frac{np}{r!(np+q)^{r+1}}\int_0^\infty{\!\!\!y^r {\rm e}^{-y} F_{n-1}\left(-\frac{y}{np+q}\right) {\rm d}y}.
\end{equation}
Using expansion of $F_n(y)$ (\ref{eq:Fdef}) in terms of central moments we get after integration the result (\ref{eq:theorem}). From Corollary~\ref{cor123} it follows that $\mu_k(n-1)$ is of the order $\mathcal{O}(n^{\lfloor \frac{k}{2}\rfloor})$ and so the $k$-th term in the expansion is $\mu_k(n-1)/(np+q)^k \sim \mathcal{O}(1/n^{\lceil \frac{k}{2}\rceil})$, showing that the expansion is indeed asymptotic in $n$.
\end{proof}
It is instructive to explicitly write out $f_r(n)$ (\ref{eq:theorem}) using first few lowest terms. We have
\begin{multline}
f_r(n)=\frac{np}{(np+q)^{r+1}} \biggl( 1+\frac{(r+2)(r+1)(n-1)pq}{2(np+q)^2}\biggl\{1 -\frac{(r+3)(q-p)}{3(np+q)}+\biggr. \biggr. \\
\biggl. \biggl. +\frac{(r+4)(r+3)(1-6pq+3npq)}{12(np+q)^2}+\cdots \biggr\} \biggr).
\label{eq:last}
\end{multline}
Because in the literature expansions in terms of inverse powers of $np$ are frequently given, instead of $(np+q)$ as here, we will write expansion in terms of $(np)$, obtained from (\ref{eq:last}) by using $1/(1+x)^k=\sum_{j=0}^\infty{(-x)^j {k-1+j \choose k-1}}$.
\begin{corollary}
$f_r(n)$ expanded in terms of inverse powers of $(np)$ is
\begin{equation}
f_r(n)=\frac{1}{(np)^r} \left(1+\frac{r(r+1)q}{2(np)} + \frac{r(r+1)(r+2)q(4+q+3rq)}{24 (np)^2}+\cdots \right).
\end{equation}
Specifically, for lowest three integer $r$'s we have
\begin{multline}
f_1(n)= \frac{1}{np} \biggl( \biggr. 1+\frac{q}{(np)}+\frac{q(1+q)}{(np)^2}+\frac{q(1+4q+q^2)}{(np)^3}+\\
 +\biggl. \frac{q(1+q)(1+10q+q^2)}{(np)^4} + \frac{q(1+26q+66q^2+26q^3+q^4)}{(np)^5}+\cdots \biggr)
\end{multline}

\begin{multline}
f_2(n)= \frac{1}{(np)^2} \biggl( \biggr. 1+\frac{3q}{(np)}+\frac{q(4+7q)}{(np)^2}+\frac{5q(1+6q+3q^2)}{(np)^3}+ \\
+\biggl. \frac{q(6+91q+146q^2+31q^3)}{(np)^4} + \cdots \biggr)
\end{multline}

\begin{multline}
f_3(n)= \frac{1}{(np)^3} \biggl( \biggr. 1+\frac{6q}{(np)}+\frac{5q(2+5q)}{(np)^2}+\frac{15q(1+8q+6q^2)}{(np)^3}+\\
 +\biggl. \frac{7q(3+58q+128q^2+43q^3)}{(np)^4} + \cdots \biggr)
\end{multline}

\end{corollary}

\section{Asymptotic Expansion for Poisson Distribution}
\label{sec2}

For Poisson distribution the whole approach is very similar to the one for binomial distribution so we will state the main Theorem right away.
\begin{theorem}
For any real $r>0$ the inverse moment $g_r$ of Poisson distribution (\ref{eq:defP}) is given by the asymptotic expansion
\begin{equation}
g_r=\frac{m}{(m+1)^{r+1}} \left\{ \sum_{k=0}^{p-1}{(-1)^k
    \frac{\tilde{\mu}_k}{(m+1)^k}{r+k \choose
      r}}+\mathcal{O}(1/m^{\lceil \frac{p}{2} \rceil}) \right\},
\label{eq:th2}
\end{equation}
where $\tilde{\mu}_k$ are central moments of Poisson distribution,
\begin{equation}
\tilde{\mu}_k=\sum_{s=0}^\infty{ {\rm e}^{-m} \frac{m^s}{s!} (s-m)^k}.
\end{equation}
\label{theoremP}
\end{theorem}
\begin{proof}
The proof essentially goes along the same steps as the one for
Theorem~\ref{theoremf}. First note that Eq.~\ref{eq:frint} can be for any
distribution written as
\begin{equation}
(1-P(s=0))\mathbb{E}(1/x^r)=\sum_{s=1}^\infty{P(s) \frac{1}{s^r}}=\frac{1}{\Gamma(r+1)}\int_0^\infty{\!\!\! x^r \frac{{\rm d}M(-x)}{{\rm d}x}{\rm d}x},
\end{equation}
where $M(x)$ is generating function of the distribution, {\em e.g.}, it is $M(x)=(q+p{\rm e}^x)^n$ for binomial and $M(x)={\rm e}^{m({\rm e}^x-1)}$ for Poisson distribution. Now we try to write derivative of generating function in terms of generating function of the central moments $F(x)$. For arbitrary distribution $F(x)$ is given by $F(x)={\rm e}^{-x\mathbb{E}(x)} M(x)$, see Lemma~\ref{lemma_F}. The
  ``trick'' is now to write ${\rm d}M/{\rm d}x$ as $a(x) F(x)$, with some
  distribution dependent $a(x)$. We can see that $a(x)$ must be given by
  logarithmic derivative of generating function, $a(x)={\rm
  e}^{x\mathbb{E}(x)} \frac{\rm d}{{\rm d}x}(\log{M(x)})$. Using $a(x)$ the
  inverse moment is then
\begin{equation}
(1-P(s=0))\mathbb{E}(1/x^r)=\frac{1}{\Gamma(r+1)}\int_0^\infty{\!\!\!x^r a(-x) F(-x){\rm d}x}.
\label{eq:genI}
\end{equation}
For Poisson distribution we have $a(x)=m{\rm e}^{x(m+1)}$, from which
  by substitution of variables $y=(m+1)x$ an expansion parameter $(m+1)$ is obtained. Furthermore, after writing series expansion for $F(-x)$ the resulting integrals $\int_0^\infty{y^p {\rm e}^{-y}{\rm d}y}$ are easily evaluated, giving final result (\ref{eq:th2}).
\end{proof}

From the proof we see that the asymptotic expansion works for distributions
having simple $a(x)$ ({\em e.g.} exponential), so that integrals
$\int_0^\infty{x^p a(x) {\rm d}x}$ can be analytically evaluated. Central moments of Poisson distribution can be calculated using the following recursive relation. 
\begin{lemma}
\label{lemmapoiss}
Central moments of Poisson distribution are given by recursive formula
\begin{equation}
\tilde{\mu}_{k+1}=m\left( \frac{{\rm d}\tilde{\mu}_k}{{\rm d}m}+k\tilde{\mu}_{k-1} \right),
\end{equation}
with the first few being
\begin{equation}
\tilde{\mu}_0=1,\quad \tilde{\mu}_1=0,\quad \tilde{\mu}_2=m,\quad
\tilde{\mu}_3=m,\quad \tilde{\mu}_4=m+3m^2.
\end{equation}
\end{lemma}

\begin{proof}
The proof goes along the same line as the one for Lemma~\ref{lemmarec}.
\end{proof}
Using Lemma~\ref{lemmapoiss} and Theorem~\ref{theoremP} we can derive first few terms in the expansion of $g_r$.
\begin{corollary}
First few term in the expansion of $g_r$ in inverse powers of $m$ are
\begin{equation}
g_r=\frac{1}{m^r}\left\{1+\frac{r(r+1)}{2m}+\frac{r(10+21r+14r^2+3r^3)}{24m^2} +\mathcal{O}(1/m^3)\right\}
\end{equation}
\end{corollary}

\section{General positive distributions}
\label{sec3}

Already in proof of Theorem~\ref{theoremP} we saw that the asymptotic expansion can be derived for general distributions. Here we state the Theorem.
\begin{theorem}
Let $M(x)$ be generating function of positive moments (about origin) for an arbitrary distribution and let $\mu_k$ denote its central moments. That is, if $\bar{x}=\mathbb{E}(x)$ then $\mu_k=\mathbb{E}((x-\bar{x})^k)$. If we denote by $\phi(x)$ the characteristic function, $\phi(x)=\log{M(x)}$, {\em i.e.} the generating function of cumulants, then the asymptotic expansion of inverse moments is in general given by
\begin{equation}
\sum_{s=1}^\infty{P(s) \frac{1}{s^r}}=\frac{1}{\bar{x}^r} \left( \sum_{k=0}^\infty \frac{(-1)^k \mu_k}{k!} \frac{1}{\bar{x}^k} \int_0^\infty{\!\!\!y^{r+k} {\rm e}^{-y}\phi'(-y/\bar{x}){\rm d}y}\right).
\label{eq:gen}
\end{equation}
The above series must be understood as an asymptotic expansion in $1/\bar{x}$.
\end{theorem}
\begin{proof}
The proof goes along the same line as the one for Theorem~\ref{theoremP}. One can start with equation~\ref{eq:genI} and after changing the integration variable to $y=\bar{x} x$ one arrives at equation~\ref{eq:gen}.
\end{proof}
Note that expansion (\ref{eq:gen}) works exactly when the asymptotic expansion makes sense, that is when {\em e.~g.} the first inverse moment is given in the leading order by the inverse of the expectation value, $1/\bar{x}$. This happens when higher central moments grow with the mean $\bar{x}$ sufficiently slowly, for instance the second central moment as $\mu_2 \propto \bar{x}^{2-\delta}$ with $\delta>0$, see also discussion in~\cite{Garcia:01}. In the case of binomial and Poisson distributions one has $\mu_k \sim \bar{x}^{\lfloor \frac{k}{2}\rfloor}$, ensuring asymptoticity of the expansion.

\section*{Acknowledgment}

The author would like to thank the Alexander von Humboldt Foundation for financial support.

\bibliographystyle{amsplain}

\bibliography{binomial}

\end{document}